\documentclass{amsart}
\usepackage{amssymb,amscd}

\newcommand\pa{\partial}
\newcommand\bS{{}^b\kern-2ptS}
\def\lra{\longrightarrow}
\def\cdc{,\ldots,}
\def\Om{\Omega}
\def\fr{\frac}
\def\ll{L^2}
\newcommand\llb{L^2_{\text{b}}}
\def\fcfm{F\in\Fs{M}}

\def\CC{\mathbb C}

\def\R{\mathbb R}
\def\ZZ{\mathbb Z}
\let\RR=\R
\def\spp{\text{ supp }}
\newcommand\clbp[1]{\mathfrak{A}_{#1}}
\newcommand\ssigma{\sigma}
\newcommand\supp{\operatorname{supp}}
\newcommand\ie{i.e.\/ }

\newcommand\bNp{\overline{N_+}}
\newcommand\bT{{}^b\kern-2ptT}
\newcommand\ci{$\mathcal{C}^{\infty}$}
\newcommand\CI{\mathcal{C}^{\infty}}
\newcommand\dCI{\dot{\mathcal{C}}^{\infty}}
\newcommand\dCIc{\dot{\mathcal{C}}_{c}^{\infty}}
\newcommand\CIc{\mathcal{C}_c^{\infty}}
\newcommand\Ind{\operatorname{Ind}}
\newcommand\Inop{\operatorname{In}}
\newcommand\Inal[2]{\Psi_{b,I}^{#1}(\bNp #2)}
\newcommand\bInal[1]{{\mathfrak A}_{#1}}
\newcommand\bnInal[1]{{\mathfrak A}_{#1}^{(-1)}}

\newcommand\Fs[1]{\mathcal F(#1)}
\newcommand\hFs[1]{{\mathcal {HF}}(#1)}
\newcommand\Diag{\operatorname{Diag}}

\newtheorem{definition}{Definition}
\newtheorem{theorem}{Theorem}
\newtheorem{proposition}{Proposition}
\newtheorem{corollary}{Corollary}
\newtheorem{lemma}{Lemma}

\begin{document}

\title[$K$-theory of b-pseudodifferential operators]
{$C^*$-algebras of {\normalsize $b$}-pseudodifferential operators
and  an $\R^k$-equivariant index theorem}

\author[Richard Melrose]{Richard Melrose$^1$}
\address{Department of Mathematics, MIT}
\email{rbm\@math.mit.edu}
\author[Victor Nistor]{Victor Nistor$^2$}
\address{Department of Mathematics, Pennsylvania State University}
\email{nistor\@math.psu.edu}

\thanks{$^1$ Partially supported by NSF grant DMS-9306389.  $^2$ Partially
supported by NSF Young Investigator Award DMS-9457859 and a Sloan research
fellowship. Preprints are available by ftp from {\bf
ftp-math-papers.mit.edu} and from {\bf http:{\scriptsize//}www-math.mit.edu%
{\scriptsize/}$\sim$rbm{\scriptsize/}rbm-home.html} or
{\bf http://www.math.psu.edu/nistor/}}

\date{\today}

\begin{abstract}
We compute $K$-theory invariants of algebras of pseudodifferential
operators on manifolds with corners and prove an equivariant index theorem
for operators invariant with respect to an action of $\R^k.$ 
We discuss the relation between our results and the $\eta$-invariant.
\end{abstract}
\maketitle

\section*{Introduction}

In this paper we determine the $K$-groups of the norm closure of the algebra,
$\Psi^0_b(M),$ of b-pseudodifferential (or totally charactertistic)
operators acting on the compact manifold with corners $M.$ In the case of a
compact manifold with boundary this class of operators was introduced in
\cite{Melrose25}, see also \cite{Melrose42} and \cite{Hormander3}. For the
general case of a compact manifold with corners it was described in
\cite{Melrose-Piazza1}. There are closely related algebras which have the
same completion, see \cite{Melrose45}.

The algebra $\Psi^0_b(M)$ can be identified with a $*$-closed subalgebra of the
bounded operators on $L^2(M)$ and the Fredholm elements can then be
characterized by the invertibility of a joint symbol, consisting of the
principal symbol, in the ordinary sense, and an `indicial operator' (as for
Fuchsian differential operators) at each boundary face, arising by freezing
the coefficients at the boundary face in question. In view of the
invariance of the index with respect to small perturbations \cite{Douglas1}
we consider (as in the case $\pa M=\emptyset$ for the Atiyah-Singer index
theorem, \cite{Atiyah-Singer1,Singer2}) the norm closure, which we denote
$\mathfrak A_M,$ of $\Psi^0_b(M).$ The closure leads to an algebra,
whose  $K$-theory is easier to  compute. Just as in the case of a manifold without
boundary, the principal symbol map has a continous extension to
$\mathfrak{A}_M,$ with values in $C(\bS^*M),$ where $\bS^*M\equiv S^*M$ as
manifolds. For any compact manifold with boundary we show that the
principal symbol map $\sigma_0:\mathfrak A_M \to C(\bS^*M)$ induces an
isomorphism in $K$-theory.

The algebra $\mathfrak A_M$ contains the compact operators
$\mathcal K.$ Denote by $Q_M=\mathfrak A_M/\mathcal K$ the quotient. If
$\pa M=\emptyset,$ then $Q_M$ is isomorphic to the algebra $C(S^*M)$ `of
symbols', in the general case we call $Q_M$ the {\em algebra of joint
symbols} since it involves both the principal symbol and the indicial
operators. In the case of a compact manifold with boundary we show that the
principal symbol map still induces an isomorphism of the $K_0$-groups,
whereas for $K_1$ the boundary, $\pa M$, through the indicial operator, contributes
an extra copy of $\mathbb Z;$ this can be attributed  to ``spectral flow''
invariants \cite{Atiyah-Patodi-Singer1}. More precisely there is a short exact
sequence 
$$
\begin{CD}
0 @>>> \mathbb Z @>>> K_1(Q_M) @>(\ssigma_{0})_*>> K_1(C_0(\bS^*M)) @>>> 0.
\end{CD}
$$
The index morphism $\Ind:K_1(Q_M) \to \mathbb Z$ provides a splitting of 
this exact sequence.

We also compare the algebraic and topological $K$-theory of the uncompleted
algebra $\Psi^{\circ}_b (M)$ and we interpret a result in \cite{Melrose46}
on the $\eta$-invariant in this setting. We conclude the paper with some
results on the index of operators on manifolds equiped with a proper
action of $\RR^k$.

In order not to have to worry about composition and boundedness of 
order zero operators, we will consider only operators whose Schwartz
kernel is compactly supported. In case our operators act on a $G$-space $X$
and are invariant with
respect to the action of $G$ (typicaly $G=\RR^k$) then we assume that
the distribution kernels have compact support in the quotient $(X\times X)/G$.

We would like to thank Is Singer for a helpful discussion.

\section{Manifolds with corners}

Below we work on a smooth manifold with corners $M.$ By definition, this
means that every point $p\in M$ has coordinate neighborhoods, which are
diffeomorphic to $[0,\infty)^k \times \R^{n-k}$ where $n$ is the dimension
of $M,$ $k=k(p)$ is the {\em codimension} of the face containing $p,$ and
$p$ corresponds to $0$ under this isomorphism. The transitions between such
coordinate neighborhoods must be smooth up to the boundary. An {\em open
face} will be a path component of the set $\partial_kM$ of all points $p$
with the same fixed $k=k(p).$ The closure, in $M,$ of an open face will be
called a boundary face, or simply a {\em face}; if it is of codimension one it
may be called specifically a boundary hypersurface. In general such a boundary
face does not have a covering by coordinate neighborhoods of the type
described above because boundary points may be identified. To avoid this
problem we demand, as part of the definition of a manifold with corners,
that the boundary hypersurfaces be embedded. More precisely this means that
we assume that for each boundary hypersurface $F$ of $M$ there is a
smooth function $\rho_F \ge 0$ on $M,$ such that
\begin{equation}
F = \{\rho_F = 0\},\ d(\rho_F)\not=0\text{ at }F.\label{BoundaryDefining}
\end{equation}
If $p \in F_0,$ a face of codimension $k,$ then exactly $k$ 
of the  functions $\rho_F$ vanish at $p.$ Denoting them $\rho_1,\ldots,\rho_k,$ then
$d\rho_1,\ldots,d\rho_k$ must be linearly independent at $p;$ it follows that the
addition of some $k$ functions (with independent differentials at $p$ on
$F_0$) gives a coordinate system near $p.$

We denote by $\Fs M$ the set of boundary faces of $M;$ in view of the assumed
existence of boundary defining functions, \eqref{BoundaryDefining}, they
are all manifolds with corners.  We can assume, without loss of generality,
that $M$ is connected, and hence that there is a unique face of codimension
$0,$ namely $M.$ For each manifold with corners we denote by $\hFs{M}$ the
set of boundary hypersurfaces $F\in\Fs{M}$ (\ie faces of codimension $1$.) 
In particular
if $F\in\Fs{M}$ then $\hFs{F}$ is the subset of $\Fs{M}$ consisting of the
$F'\in\Fs{M}$ which are contained in $F$ and with codimension equal to the
codimension of $F$ plus one.

We now introduce some constructions and notation to be used in the sequel.

We make a choice of functions $\rho_F$ as in \eqref{BoundaryDefining} and 
fix a metric $h$ that locally at any point $p$ has the form 
$h = (dx_1)^2 +\dots+ (dx_k)^2 + h_0 (y_1,\ldots,y_{n-k})$ if 
$x_1,\ldots,x_k,$  $y_1,\ldots,y_{n-k}$ are local coordinates at $p,$ and 
$x_1=\rho_{F_1},\ldots,x_k=\rho_{F_k}$ are defining functions as 
in our assumption. The existence
of such a metric is shown in \cite{Hassell-Mazzeo-Melrose2}, for example.

The choice of the functions $x_{H}=\rho_H$ for all $H\in\hFs{M}$ establishes a
trivialization $N F \simeq F \times \R^k$ of the normal bundle to each
boundary face.  We denote by $N_+ F \subset N F$ the closed set of normal
vectors that point into the manifold $M.$ They are exactly those vectors
which have non-negative $x$-components. Now the group 
$(0, \infty)^k = \R^{*k}_+$ acts on
$\ll (N_+ F, \Om^{\fr12}_b)$ by dilations
\begin{multline}
\label{action} 
\lambda_\epsilon(u)(x_1,\ldots,x_k, y_1,\ldots,y_{n-k})\\
=        u (\epsilon_1^{-1} x_1,\dots, \epsilon^{-1}_k x_k,
y_1,\ldots,y_{n-k}),\; 
	\epsilon=(\epsilon_1,\dots,\epsilon_k).
\end{multline}
This action is independent of the choice of defining functions.

The exponential map associated to the
Levi-Civita connection gives a diffeomorphism from a neighborhood
$V_F$ of the zero section in $N_+ F$ to an open neighborhood of $F$ in
$M$
$$
	\Phi_F = \exp: V_F \longrightarrow M,\ V_F\subset N_+  F.
$$
Due to the particular choice of the metric $h,$ $\Phi_F$ is a
diffeomorphism of manifolds with corners which maps the zero section of
$NF$ onto $F.$

Let $\varphi_F$ be a smooth function on $M,$ $0 \le \varphi_F \le 1,$ 
supported inside $\Phi_F (V_F),$ and such that 
$\varphi_F=1$ in a neighborhood of $F.$
We define 
\begin{equation}
	L_F: \ll (N_+ F, \Om^{\fr12}_b) = \ll ( N_+ F, 
	 \fr{dx_1\ldots dx_k}{x_1\ldots x_k}
	dy_1\ldots dy_{n-k}) \longrightarrow \ll (M, \Om^{\fr12}_b)
\label{L2insertion}
\end{equation}
by the formula $L_F(u) = \varphi_F (u \circ \Phi_F^{-1}),$ 
well defined since $\spp \varphi_F \subset \Phi_F(V_F).$

\section{The algebra of b-pseudodifferential operators}

A definition of the algebra of b-pseudodifferential operators is recalled
in the appendix. The definition given there starts from an explicit
description of the algebra $\Psi_b^*(M)$ for the special, model, case of
$M=[-1,1]^n.$ The general case is then obtained by localization and
$\Psi_b^*(M)$ consists of operators from $\CIc(M)$ to $\CI(M).$ If $M$ is a
manifold without boundary this definition reduces to that of
$1$-step polyhomogeneous (i.e. classical) pseudodifferential operators in
the usual sense.

As in the boundaryless case, the principal invariant of a
pseudodifferential operator is its {\em principal symbol}, it is a function
on the b-cotangent space. This bundle, denoted $\bT^*M,$ is naturally
defined over any manifold with corners. Over the interior it is canonically
identified with $T^*M$ but at a boundary point, $p,$ its fibre is the space of
equivalence class of differentials 
\begin{equation}
\sum\limits_{p\in H}a_H\frac{d\rho_H}{\rho_H}+d\phi,\ a_H\in\RR,\ 
H \in \hFs M,\ \phi\in\CI(M)
\label{SectionofbT*M}\end{equation}
modulo the space of smooth differentials which, after being pulled back to 
$F=\bigcap\limits_{p\in H}H,$ vanish at $p.$ It can be defined more
naturally as the dual bundle to the bundle, $\bT M,$ with sections
consisting precisely of the space of smooth vector fields on $M$ tangent to
all boundary hypersurfaces.

Let $\bS^*M$ be the quotient of $\bT^*M\setminus0$ by the fibre action of
$(0,\infty)$ and let $P^m$ be the bundle over $\bS^*M$ with sections which
are homogeneous functions of degree $m$ on $\bT^*M.$

\begin{proposition}\label{PseudoSymbolMap} There is a natural short exact
sequence 
\begin{equation}
0\longrightarrow \Psi_b^{m-1}(M)\hookrightarrow
\Psi_b^m(M)\overset{\sigma_m}{\longrightarrow}
\CI(\bS^*M;P^m)\longrightarrow 0,
\label{bPseudoSymbol}
\end{equation}
which is multiplicative if $M$ is compact, where $\sigma_m(A)$ is
determined by `oscillatory testing' in the sense that if $\psi\in\CIc(M),$
$\phi\in\CI(M)$ is real valued and $a_H\in\RR$ are such that the
corresponding section $\alpha$ of $\bT^*M$ given by \eqref{SectionofbT*M}
is non-vanishing over the support of $\psi$ then
\begin{equation}
\sigma_m(\alpha)\psi=
\lim\limits_{\lambda\to\infty}\lambda^{-m}\prod\limits_{p \in H}\rho_H^{i\lambda
a_H}e^{i\lambda \phi}A\big(\prod\limits_{p \in H}\rho_H^{-i\lambda a_H}e^{-i\lambda
\phi}\psi\big).
\end{equation}
\end{proposition}

Consider the projective compactification of the closed half-line as in
\eqref{ProjectiveCompactification}
\begin{equation}
[0,\infty)\ni s\longrightarrow \frac{s-1}{s+1}\in [-1,1].
\end{equation}
The multiplicative action of $(0,\infty)$ on $[0,\infty)$ lifts to be
smooth on $[-1,1]$ so the $k$-fold application of this compactification
embeds the inward pointing normal bundle $N_+F=[0,\infty)^k\times F$ to any
boundary face of a manifold with corners, into $\bNp F=[-1,1]^k\times F,$
with the \ci\ structure on the compactification independent of the choice
of boundary defining functions used to produce the trivialization; the
action of $\R^{*k}_+$ lifts to be smooth on $\bNp F.$

\begin{definition} If $M$ is a compact manifold with corners then the indicial
algebra $\Inal*F$ corresponding to a boundary face $F \subset M$
is the algebra consisting of those b-pseudodifferential operators $T$ on
$\bNp F$ which are invariant under the $\R^{*k}_+$ action (see the remark
on invariance following Definition~\ref{BasicDefinition} in  the
appendix).
\end{definition}

The operators $T \in \Inal*F$ of order at most $m$ form a subspace denoted
$\Inal mF,$ so $\Inal*F=\bigcup_m\Inal mF.$  Let $L_F$ and the action
$\lambda$ be as in \eqref{action} and \eqref{L2insertion}.
Recall the following result \cite{Melrose42}

\begin{theorem} \label{Theorem.IndicialMaps} For any boundary face $\fcfm$
there is an onto  morphism
$$
\Inop_{F,M}:\Psi^0_b (M) \lra\Inal0{F},\
$$
independent of any choices and uniquely determined by the property
$$
\Inop_{F,M}(T)u=\lim\limits_{\epsilon_i\to 0}(\lambda_{\epsilon^{-1}} L^*_F T L_F
\lambda_\epsilon )u
$$
for any $u \in\CIc(N_+F),$ $\epsilon=(\epsilon_1,\ldots,\epsilon_k),$ $k$
being the codimension of $F.$
\end{theorem}

\begin{proof} Suppose that $T\in \Psi^0_b (M).$ As discussed in the
Appendix, $T$ is locally of the form \eqref{SimpleAction}. If
$x_1,\ldots,x_k$ are defining functions for the face to which
$p$ belongs and $y_1,\dots,y_{n-k}$ are additional local coordinates then
the defining formula  \eqref{SimpleAction} reduces to
\begin{multline}
Tu(x,y)=\int\limits_0^\infty\dots\int\limits_0^\infty\int\limits_{\RR^{n-k}}
T(x_1,\dots,x_k,x'_1,\dots,x'_k,y_1,\ldots,y_{n-k},y'_1,\dots,y'_{n-k})\\
u(x'_1x_1,\dots,x'_kx_k,y'_1,\dots,y'_{n-k})
\frac{dx'_1}{x'_1}\dots\frac{dx'_k}{x'_k}dy'_1\dots dy'_{n-k}
\end{multline}
where now $T(s,x,y,y')$ is conormal at $x'_i=1,$ $y=y'$ or smooth as the
localizing functions are in the same or different coordinate patches; it is
still rapidly decreasing as $x'_i\to0$ or $\infty$ and now has compact
support in $y,y'.$

Since the computation is local we can assume that  $T=L_F^*T L_F$ and then
\begin{multline}
(\lambda_{\epsilon^{-1}}T \lambda_\epsilon)(u)(x,y)\\
=\int\limits_0^\infty\dots\int\limits_0^\infty\int\limits_{\RR^{n-k}}
T(\epsilon_1x_1,\dots,\epsilon_kx_k,
x'_1,\dots,x'_k,y_1,\ldots,y_{n-k},y'_1,\dots,y'_{n-k})\\
u(x'_1x_1,\dots,x'_kx_k,y'_1,\dots,y'_{n-k})
\frac{dx'_1}{x'_1}\dots\frac{dx'_k}{x'_k}dy'_1\dots dy'_{n-k}
\end{multline}
(after a dilation in the $x'$-variables).
This shows immedately that the limit from the statement
exists  as $\epsilon\to0$and that the the localized indicial operator 
is given by
\begin{multline}
\label{eq10}
\Inop_{F,M}(T)u=
\int\limits_0^\infty\dots\int\limits_0^\infty\int\limits_{\RR^{n-k}}
T(0,\dots,0,x'_1,\dots,x'_k,y_1,\ldots, y_{n-k},y'_1,\dots,y'_{n-k})\\
u(x'_1x_1,\dots,x'_kx_k,y'_1,\dots,y'_{n-k})
\frac{dx'_1}{x'_1}\dots\frac{dx'_k}{x'_k}dy'_1\dots dy'_{n-k}
\end{multline}
for any $u\in\dCI(N_+F).$ Conversely, if $T$ is as in \eqref{SimpleAction}
and does not depend on $x$, the above formula defines an invariant operator
in $\bInal0F$. This shows that $\Inop_{F,M}$ is onto. 
\end{proof}

We make the observation that the form of the operators $T \in\Psi_b^0(M)$
ensures that $Tu\vert_{F}$ depends only on $u\vert_F,$ so $T$ restricts 
to an operator $T\vert_F$ on $F.$ The indicial operator $\Inop_{F,M}(T)$
captures not only this restriction but all restrictions $(x^{i\lambda}T
x^{-i\lambda})\vert_F$ for $\lambda\in\CC.$ This is the original definition
of the indicial operators \cite{Melrose42}. 
The equivalent definition for order zero operators contained in the previous theorem
has the advantage that it readily extends to the completion in norm.

\section{The norm closure, $\mathfrak{A}_M$}

We now proceed to discuss the algebra, $\mathfrak A_M,$ obtained by taking
the norm closure of $\Psi^0_b(M)$ as an algebra of bounded operators on
$\llb(M).$ As in the case of a compact manifold without boundary the symbol
map extends by continuity.

The elements of $\Inal0F$ define bounded operators on $\llb (N_+
F);$ the norm closure of this subalgebra of $\Psi_b^0(\bNp F)$ will be
denoted by $\bInal{F,M}.$   
The closure in norm of $\Inal{-1}F$ will be denoted $\bnInal{F,M}$.
Thus $\mathfrak A_M=\mathfrak A_{M,M}$.
When $M$ is understood we will omitt the last $M$ and write
$\bInal{F}=\bInal{F,M}$ and $\bnInal{F}=\bnInal{F,M}$.

The following proposition extends the indicial morphisms to the norm
closed algebras introduced above. 

\begin{proposition} 
\label{ExtIndCont}
The indicial morphisms extend to maps 
$
\Inop_{F,M}: \clbp{M}\lra \bInal{F}
$
for any boundary face $F$ of $M$ and by iteration
to  morphisms 
$$
\Inop_{F',F}:\bInal{F}\lra\bInal{F'}
$$ 
for any boundary face $F' \subset F.$  These morphisms decrease the norm, 
are uniquely determined and satisfy 
$
\Inop_{F'',F'} \circ \Inop_{F',F} = \Inop_{F'',F}
$ 
for all boundary faces $F'' \subset F' \subset F.$  
\end{proposition}

\begin{proof} It follows from the definition 
of the indicial morphisms $\Inop_{F,M}$ given in 
Theorem~\ref{Theorem.IndicialMaps} that they satisfy
$$
\| \Inop_{F,M}(T)u \| = \| \lim\limits_{\epsilon_i\to 0}
(\lambda_{\epsilon^{-1}} L^*_F T L_F \lambda_\epsilon )u \| \leq \|T\| \,\|u\|
$$
and hence $\| \Inop_{F,M}(T)\| \leq \|T\|$. This show that $\Inop_{F,M}$
extends by continuity to the norm closure.

We now define the morphisms $\Inop_{F',F}$ for $F' \subset F$. Denote by
$\epsilon = (\epsilon',\epsilon'')$ corresponding to local coordinates
$(x',x'',y)$ where $x'$ are defining functions for $F$ and $x=(x',x'')$
are defining functions for $F'$. The relation 
\begin{multline}
\Inop_{F',M}(T)u=\lim\limits_{\epsilon_i\to 0}
(\lambda_{\epsilon^{-1}} L^*_{F'} T L_{F'}
\lambda_\epsilon )u\\
=\lim\limits_{\epsilon_i''\to 0}\big(\lambda_{\epsilon''}^{-1}
\lim\limits_{\epsilon_i'\to 0}\big(\lambda_{\epsilon'}^{-1} 
T \lambda_{\epsilon'} \big ) \lambda_{\epsilon''}\big )u
=\lim\limits_{\epsilon_i''\to 0}\left(\lambda_{\epsilon''}^{-1}
\Inop_{F,M}(T) \lambda_{\epsilon''}\right ) u
\end{multline}
if $u$ has very small support shows that $\Inop_{F',M}(T)$
depends only on $\Inop_{F,M}(T)$ and that 
$\|\Inop_{F',M}(T)\| \leq \|\Inop_{F,M}(T)\|$. Since $\Inop_{F,M}$
is onto, Theorem \ref{Theorem.IndicialMaps}, there exists a unique
contraction $\Inop_{F',F}$ satisfying 
\begin{equation}
\label{DefInd}
\Inop_{F',F} \circ \Inop_{F,M} = \Inop_{F',M}
\end{equation}
The compatibility relation 
$\Inop_{F'',F'} \circ \Inop_{F',F} = \Inop_{F'',F}$ follows from
uniqueness in the definition \eqref{DefInd}. This proves the proposition.
\end{proof}

For a locally compact space $X$ we shall denote by $C_0(X)$ the algebra of
those continuous functions on $X$ which vanish at infinity. It is the norm
closure of the algebra $C_c(X)$ of compactly supported functions on $X.$
If $X$ is a smooth manifold the set of compactly supported smoot functions
will be denoted by $\CIc(X)$ and is also dense in $C_0(X)$.

We now construct a cross-section for  $\Inop_{F'F}.$

\begin{proposition} There exist linear maps 
$\lambda_{F ,F'}:\bInal{F'} \lra \bInal{F}$ 
defined for all pairs $F', F$ where $F'$ is a face of $F$ of codimension 1, 
with the following properties

(i)  $\lambda_{F, F'} (\Inal0{F'}) \subset \Inal0{F},$

(ii)  $\Inop_{F',F} \circ \lambda_{F,F'} (T) = T$ for all $T \in\bInal{F'},$

(iii) $\|\lambda_{F ,F'} (T) \| \le \| T \|,$

(iv) Let $F''$ be an other face of $F,$ not contained in $F',$ then
$$
	\Inop_{F'', F} \circ \lambda_{F, F'} (T) = 
\cases 
	\lambda_{F'', F' \cap F''} \circ \Inop_{F'\cap F'',F'} (T) 
 	& \qquad \text{if} \;\;F' \cap F''\not = \varnothing \\
	0 & \qquad \text{otherwise}.
\endcases
$$
\end{proposition}

\begin{proof} \ We observe  that $F'\cap F''\not =\varnothing$ 
exactly when $F'\cap F''$ is a face.
We first define $\lambda_{F_0, M}$ for every boundary hypersurface  
$F_0$ of $M$ by the formula $\lambda_{M, F_0} (T) = L_{F_0} T L^*_{F_0}.$  
In local coordinates $\lambda_{M, F_0}(T)$ is given by the formula
$\lambda_{M,F_0}(T)(t,s)=\phi_{F_0}(t)T(t,s)\phi_{F_0}(s)$, 
where $\phi_{F_0}$ is as in equation \eqref{L2insertion}.

For an arbitrary pair $(F, F')$ as in the statement there exists a unique 
maximal face (\ie boundary hypersurface) $F_0$ such that $F' = F \cap F_0.$
Using local coordinates as above we see that by construction the indicial operator 
$\Inop_{F, M} \circ \lambda_{M, F_0} (T)$ depends only on 
$\Inop_{F', F_0} (T).$  Using this fact we see that we can define the 
linear section $\lambda_{F, F'}$ by the formula 
\begin{equation}
\label{DefSect}
	\lambda_{F, F'}\circ \Inop_{F', F_0} (T) =
        \Inop_{F, M} \circ \lambda_{M, F_0} (T)
\end{equation}

Consider now three faces $F, F'$ and $F'',$ $F',F'' \subset F,$ as in the 
statementof the proposition.
Let $F_0$ be the unique maximal face of $M$  such that 
$F' = F \cap F_0,$ as above.  Then we have for all $T \in\bInal{F_0}$

\noindent $\Inop_{F'',F} \circ \lambda_{F, F'} \circ \Inop_{F', F_0} (T)  = 
\Inop_{F'', F} \circ \Inop_{F, M} \circ \lambda_{M, F_0} (T) =
\Inop_{F'', M} \circ \lambda_{M, F_0} (T) = $
$$=
\cases 
	0  & \qquad \text{  if  }  F'' \cap F_0 = \varnothing \\
	%
	\lambda_{F'', F'' \cap F_0} \circ \Inop_{F'' \cap F_0, F'} 
	\circ \Inop_{F', F_0} (T) & \qquad \text{ otherwise }
\endcases
$$
where we have used  the definition \eqref{DefSect}, the properties of the 
indicial  morphisms proved in the previous proposition and the fact that 
$F'' \cap F_0 = (F'' \cap F) \cap F_0 = F'' \cap F'.$
This is enough to conclude the proof in the case  $F'' \cap F' = \varnothing.$  

If $F'' \cap F' \neq \varnothing$ 
we also need to observe  
that $F'' \cap F'$ is of codimension 1 in $F''$ and hence 
$\lambda_{F'', F' \cap F''}$ is defined.
\end{proof}

{\em Observation}. It is possible to by putting some extra conditions on the 
functions $\phi_{F}$ to  define $\lambda_{F',F}$ satisfying
$\Inop_{F',F} \circ \lambda_{F',F}=id$ and 
$\lambda_{F'',F'}\circ\lambda_{F',F}=\lambda_{F'',F}$.

\begin{corollary} \label{Cor1}\ Suppose the operators 
$T_{F'} \in \bInal{F'}$, respectively 
$T_{F'} \in \Inal0{F'}$, $F' \in \hFs F,$ satisfy the 
compatibility condition
$$
\Inop_{F' \cap F'', F'} (T_{F'}) = \Inop_{F' \cap F'', F''} (T_{F''}) 
$$
for all pairs $F', F''.$   Then we can find $T \in \bInal{F}$, 
respectively  $T \in \Inal0{F}$, such that 
$T_{F'} = \Inop_{F', F} (T)$ and 
$\| T \| \le C \max \; \|T_{F'}\|,$
where the constant $C >0$ depends only on the face $F.$
\end{corollary}

\begin{proof} Let $\hFs F = \{F_1, F_2\cdc F_m\}$ and define
$
	T_1 = \lambda_{F, F_1} (T_{F_1}) \in \bInal{F}
$
(respectively $T_1\in \Inal0{F}$) and 
$
	T_{l+1} = T_l + \lambda_{F,F_{l+1}} (T_{F_{l+1}} - 
	\Inop_{F_{l+1},F} (T_l)).
$
We will prove by induction on $l$ 
that $\Inop_{F_j, F} (T_l) = T_{F_j}$ for all indices 
$j \le l.$  Indeed for $l = 1$ this is the basic property of the 
sections $\lambda.$  We now prove the statement for $l +1$ and $j \le l.$

If $F_j \cap F_{l +1}=\varnothing$ we have 
$\Inop_{F_j,F}\circ \lambda_{F,F_{l +1}}=0.$
If $F_j \cap F_{l +1}\not =\varnothing$ we have 
$
	\Inop_{F_j,F}\circ \lambda_{F,F_{l +1}}
        =\lambda_{F_j, F_j \cap F_{l+1}} 
	\Inop_{F_{j} \cap F_{l+1}, F_{l+1}}.$
Using the induction hypothesis for $l$ and the compatibility relation from the
assumptions of the Corollary we have 
$$
\aligned
\Inop_{F_{j} \cap F_{l+1}, F_{l+1}} 
	(T_{F_{l+1}} - \Inop_{F_{l+1}, F} (T_l)) & =
	\Inop_{F_{j} \cap F_{l+1}, F_{l+1}} (T_{F_{l+1}}) -  
	\Inop_{F_{j} \cap F_{l+1}, F} (T_{l}) \\
	= & \Inop_{F_{j} \cap F_{l+1}, F_{l+1}} (T_{F_{l+1}}) - 
	\Inop_{F_j \cap F_{l+1}, F_j} \circ 
	\Inop_{F_j, F} (T_{l}) \\
	= & \Inop_{F_{j} \cap F_{l+1}, F_{l+1}} (T_{F_{l+1}}) -    
	\Inop_{F_{j} \cap F_{l+1}, F_j} (T_{F_j})  =  0
\endaligned
$$
Using this computation we obtain in both cases 
$
\Inop_{F_j, F} (T_{l+1}) = \Inop_{F_{j},F} (T_{l})=T_{F_j}.$
Finally, for $j = l+1,$\ $\Inop_{F_j,F}\circ \lambda_{F,F_{l +1}}=id$ and
we obtain:
$$ 
\Inop_{F_{l+1, F}} (T_{l+1}) = 
\Inop_{F_{l+1, F}} (T_{l}) - T_{F_{l+1}} -  
\Inop_{F_{l+1, F}} (T_{l}) = T_{F_{l+1}}\,.
$$  
From construction we see that our corollary is satisfied if we take $C=3^{m-1}$.
\end{proof}

\begin{lemma} \label{lemma1}
Let $X$ be a smooth manifold. If $f$ is a smooth compactly 
supported function on $S^*X$ then we can find $T \in \Psi^0(X)$
with compactly supported Schwartz kernel such that $\sigma_0(T)=f$ and $\|T\|=\|f\|$.
If $T_n$ is a sequence of operators $T_n \in \Psi^0(X)$ whose Schwartz
kernels are compactly supported and satisfy $\|T_n-T\| \to 0$ and 
$\|\sigma_0(T_n)\|\to 0$ then we can choose $T_n'  \in \Psi^{-1}(X)$
such that $T_n'$ converges to $T$  and the Schwartz kernels of $T_n'$ have 
compact support.
\end{lemma}

\begin{proof} Assume first that $X$ is compact. The first statement
is then a classical fact. It follows from example from  a theorem of Epstein
and Melrose \cite{Epstein-Melrose-Mendoza1} identifying
pseudodifferential operators with Toeplitz operators associated to the 
Dolbeaux operator on $T^*X$ on a small tubular neighborhood of $X$. 
In case $X$ is not compact but $f$ has
compact support we can still use the Epstein-Melrose theorem
for the double of $X_0$ where $X_0$ is a compact manifold with boundary
containing (the projection of) the support of $f$ in its interior.

The second statement is an immediate consequence of the first statement.
Indeed, choose using the first part of the lemma
an operator $R_n$ satisfying $\sigma_0(R_n)= \sigma_0(T_n)$
and $\|R_n\| =\| \sigma_0(T_n)\|\to 0$. Then define $T_n'=T_n-R_n$. 
\end{proof}

\begin{definition} The principal symbol morphisms
\begin{equation}
\label{bExtPseudoSymbol}
	\sigma_{F}=\ssigma_{0F}: \Inal0F \lra C_c^\infty (\bS^* M \vert_F)
\end{equation}
are defined using the identification of $\R_+^{*k}$-invariant  
functions on $\bS^* N_+F$ 
with their value at the zero section:
$C(\bS^* N_+F)^{\R_+^{*k}} \simeq C(\bS^* N_+F \vert_F) \simeq C(\bS^* M \vert_F).$
\end{definition}

In the above definition all the $b$-cotangent bundles are
(non-canonicaly) isomorphic to the ordinary cotangent bundles
$S^*N_+F$, $S^* M$. 

The normed closed algebras that we consider, $\bInal{M}$, $\bInal{F},$
$\bnInal{M}$ and $\bInal{F}$ are all closures of algebras acting
on a Hilbert space. Moreover since the algebras we start with,
$\Psi_b^0(M)$, $\Inal{0}F$, $\Psi_b^{-1}(M)$, $\Inal{-1}F$ are
closed under $*$ ($T^*$ is the Hilbert space adjoint of $T$, 
\cite{Douglas1}), the algebras we obtain are also closed in norm. Norm
closed algebras acting on  a Hilbert  space and closed under involution
are called $C^*$-algebras. This is actually a theorem due to Gelfand
and Naimark. The actual definition of a $C^*$-algebra is that of
a Banach algebra equiped with an involution $*$ such that $\|T^*T\|=\|T\|^2$
for any $T$ in the algebra. Bellow we will use the fact that any 
{\em algebraic} morphism of $C^*$-algebra is continuous and has
closed range \cite{Diximier1}.

Recall that all our opertors have compactly
supported kernels modulo the action of $\RR^k$. 

\begin{lemma}\label{SymbolMap} The symbol map \eqref{bExtPseudoSymbol}
extends by continuity to a surjective map 
$$
\ssigma_{0F}: \bInal{F} \lra C_0 (\bS^* M \vert_F).
$$
\end{lemma}

\begin{proof} The principal symbol morphism $\ssigma_{0,F}$
is a $*$-morphism, \ie it satisfies 
$$
\sigma_{0,F}(T^*)=\overline{\sigma_{0,F}(T)}, 
$$
and consequently its range is closed \cite{Diximier1}.
Moreover its  range contains $C_c^\infty(\bS^*M\vert_{F})$ which is
dense so $\sigma_{0,F}$ is onto. It follows from the oscillatory testing property
of the principal symbol map, Proposition \ref{PseudoSymbolMap}, that
$\|\sigma_0(T)\| \leq \|T\|$ for all $T \in \Psi_b^0(M)$. 
\end{proof}

Fix a face $F$ of $M$,  possibly $M$ itself.
We are going to use the previous results, especially
the last  corollary to study the indicial algebras $\Inal{F}.$
Denote by $\sigma_F^{tot}$ the 'joint symbol morphism' 

\begin{equation}
\label{CompF}
	\sigma_F^{tot}:\oplus_{F'}\Inop_{F', F} \oplus \ssigma_{0F}: 
        \bInal{F} \lra 
	\oplus_{F'}\bInal{F'} \oplus C_0 (\bS^*M\big|_F)\,,\;\; F' \in \hFs  F.
\end{equation}
Corollary \ref{Cor1} identifies the range of this morphism 
as the set of operators satisfying the `obvious' compatibility conditions:
\begin{multline*}
\Im(\sigma_F^{tot})=\{(T_{F'},f) \in 
\oplus_{F'}\bInal{F'} \oplus C_0 (\bS^*M\big|_F), \\
\Inop_{F'',F}(T_{F'})=\Inop_{F',F}(T_{F''}), \text{ and } 
f\vert_{F'}=\sigma_{0F}(T_{F'})\}.
\end{multline*}

Denote by ${\mathfrak K}_F\subset \bInal{F}$ the kernel of the joint symbol
morphism $\sigma_F^{tot}$
\begin{equation}
\label{FrakK}
{\mathfrak K}_F=\{T \in \bInal{F}, \Inop_{F',F}(T)=\sigma_{0F}(T)=0 
\;,\,\forall F' \in {\mathcal {HF}} \} .
\end{equation}

One of the basic properties of the b-calculus is that the map 
\begin{equation}
\label{automorphism}
	m_{\lambda}  (T) = \rho_{F'}^{\lambda} T \rho_{F'}^{-\lambda}\, , \;\; 
\lambda \in \mathbb C
\end{equation}
defines a grading preserving automorphism of the algebras $\Inal0{F}$
\cite{Melrose42}, for any defining function $\rho_{F'}$.  Also 
we have exact sequences
\begin{equation}
\label{eq-x}
	0 \lra \rho_F \Inal{m}{F} \lra  \Inal{m}{F}
	\overset{\Inop_{F', F}}{\!\!--\!\!-\!\!\!\!\lra} \Inal{m}{F'}\lra 0
        \, , \;\; F' \in \hFs F.
\end{equation}

Let $(\psi_m)_{m\ge 1},$ such that $0\leq \psi_m \leq 1$,  
be an increasing sequence of compactly supported
smooth functions in the interior of  $F$, $\lim \psi_m=1$ 
everywhere. Using any trivialization of $NF$ we obtain, by projection,
a similar exhausting sequence of functions supported
inside $NF,$ denoted also $\psi_m.$ This sequence is chosen such that
the projection $NF \to F$ maps the support of $\psi_m$
to a compact subset of $F$.

\begin{lemma} \ \label{lemma3} 
Define  an exhausting sequence,  $\psi_m$, of $NF$ as
above, $\psi_m=\psi_m'\circ p$, $p:NF \to F$, $\psi_m\nearrow 1$. 
Then the union  $\cup_m \psi_m \Inal{-1}{F}\psi_m$ is dense in 
${\mathfrak K}_F$ in the norm topology of $\bInal{F}.$ 
\end{lemma}

\begin{proof} Let $T \in  {\mathfrak K}_F$. 
We need  to show that we can choose the sequence $T_n \in \cup_m \psi_m 
\Inal{-1}{F}\psi_m$, $T_n \to T$.
Choose first a sequence $T_n \in \Inal0F$, $T_n \to T$.
The continuity of the indicial maps and of the pricipal symbol morphism
(Proposition \ref{ExtIndCont} and Lemma \ref{SymbolMap}) show that
$\Inop_{F', F} (T_n) \lra \Inop_{F', F} (T) = 0$ and 
$\ssigma_{0F} (T_n) \to \ssigma_{0F} (T) =  0,$.

We infer from the  corollary  \ref{Cor1} applied to $T_{F'}=\Inop_{F',F}(T_n)$
that we can find a sequence $R_n \in \Inal0F,$ 
$\Inop_{F', F} (R_n)=\Inop_{F', F} (T_n)$, $
\| R_n \| \le C \max \; \|\Inop_{F',F}(T_n)\|$. It follows that
$\|R_n\| \to 0$. Replacing
$T_n$ by $T_n - R_n$ we see that we can assume that 
$\Inop_{F',F} (T_n)=0.$ Denote by $f = \prod \rho_{F'}$ the product of all
defining functions for all boundary hypersurfaces $F'$ of $F.$ 
Using the exact sequence \ref{eq-x} we obtain that we can factor 
$T_n = f T'_n$ for some $T'_n \in \Inal{0}{F}.$  Using the invariance 
of the b-calculus with respect to the inner automorphisms 
$m_{\lambda}$ defined above, we see that
$
	T_n=f^{\fr12} T''_n  f^{\fr12}
$
for some $T''_n \in \Inal{0}F.$  The operators 
$P_{n,m}=\psi_m T_n ''\psi_m$ satisfy
$\ssigma_{0F}(P_{n,m}) \to 0$ and are usual pseudodifferential operators
on an open manifold (not b-pseudodifferential ones). 
Using  lemma \ref{lemma1} we can find
pseudodifferential operators $R_{n,m}\to 0,$ with the same principal symbol
as $P_{n,m}.$ Using  
$f^{\fr12} = \displaystyle{\lim_{n \to \infty}} \; \psi_m f^{\fr12}$ 
in the uniform norm, we have:
$$
\displaystyle{\lim_{m \to \infty}}\ \  T_n 
- f^{\fr12}(\psi_mT_n''\psi_m - R_{n,m})f^{\fr12}=
\displaystyle{\lim_{m \to \infty}}\ \  T_n 
- f^{\fr12}\psi_m(T_n'' - R_{n,m})\psi_m f^{\fr12}=0\;\,\forall n.
$$
Moreover $\ssigma_{0,F}(f^{\fr12}(\psi_mT_n''\psi_m - R_{n,m})f^{\fr12})=0.$
This shows that we can choose a sequence $m_n$ and replace 
$T_n$ by $f^{\fr12}\psi_{m_n}(T_n'' - R_{n,m_n})\psi_{m_n}f^{\fr12}$
to  obtain a sequence with the desired properties.
\end{proof}

Denote by ${\mathcal K}_F$ the algebra of compact operators on
$L^2(F)$, and by $C_0 (\R^k, {\mathcal K}_F)$ the algebra of 
${\mathcal K}_F$-valued 
functions on $\R^k,$ vanishing at infinity.

\begin{proposition} \label{Prop5} The Fourier transform in the fiber direction 
establishes an isomorphism ${\mathfrak K}_F \simeq C_0(\R^k, {\mathcal K}_F)$
for each face of codimension $k.$
\end{proposition}

If $F$ has dimension zero, \ie if $F$ reduces to a point, then 
${\mathcal K}_F=\CC$.

\begin{proof} Let $T \in \psi_m \Inal{*}F \psi_m.$  The action 
of $T$ and its adjoint $T^*$ are trivial at all boundaries of 
$N_+ F$ except $F$ itself, so we can assume 
$F = \R^{n-k},$ $N_+ F = \R^k_+ \times \R^{n-k}.$ Denote by
${\mathcal K}={\mathcal K}_F$ the algebra of compact operators on $\RR^{n-k}$.

A b-pseudodifferential operator $T$ on $N_+ F$ is given locally by
\eqref{SimpleAction} which in local coordinates 
$(x_1,\dots,x_k,y_1,\ldots,y_{n-k}) \in [0,1)^k \times \RR^{n-k}$ reads

\begin{multline}
Tu(x_1,\dots,x_k,y_1,\ldots,y_{n-k})=\\
\int\limits_0^\infty\dots\int\limits_0^\infty
A(x_1,\dots,x_k,y_1,\ldots,y_{n-k},x'_1,\dots,x'_k,y'_1,\ldots,y'_{n-k})\\
u(x_1x'_1,\dots,x_kx'_k,y'_1,\ldots,y'_{n-k})
\frac{dx'_1}{x'_1}\dots\frac{dx'_k}{x'_k}dy'_{1}\ldots dy'_{n-k}
\label{SimpleAction2}
\end{multline}
Here $A$ is a distribution as on  $[0,1)^k \times \RR^{n-k} \times 
[0,1)^k \times \RR^{n-k}$ with compact support,
%
%
with singular support contained in the set $D_0=\{x'_i=1\}\cap 
\{y_i=y'_i\}$ and with conormal
singularities at $D_0$. This means that the functions
$$
(x,y,\xi',\zeta) \to 
\int\limits_{\RR^k} \int\limits_{\RR^{n-k}} 
e^{\imath\xi' \cdot x' +\imath \zeta\cdot(y-y')} A(x,y,x',y') dx'dy'
$$
is a smooth function which has an assymptotic expansion in 
homogeneous functions of $(\xi',\zeta)$.

If the operator $T$ is invariant then we drop the condition that it has compact 
support and replace it by the condition that $\supp(A)\subset [0,1)^k \times L 
\times [0,1)^k \times L$ where $L$ is a compact subset of $\RR^{n-k}$. Then
the invariance condition is that $A$ does not depend on 
$x_1,\ldots,x_k$, so the above formula can be writen in this case as

\begin{multline}
Tu(x_1,\dots,x_k,y_1,\ldots,y_{n-k})=\\
\int\limits_0^\infty\dots\int\limits_0^\infty
A(y_1,\ldots,y_{n-k},x'_1,\dots,x'_k,y'_1,\ldots,y'_{n-k})\\
u(x_1x'_1,\dots,x_kx'_k,y'_1,\ldots,y'_{n-k})
\frac{dx'_1}{x'_1}\dots\frac{dx'_k}{x'_k}dy'_{1}\ldots dy'_{n-k}
\label{SimpleActionInv}
\end{multline}

Consider the unitary transformation
$$
	W: \ll (NF) \lra \ll (N_+ F, \Omega_b^{\fr12})
$$
given by $W f (x_1,\ldots, x_k, y) = f (\ln  x_1 \cdc \ln x_n, y) 
(x_1\ldots x_k)^{-\fr{1}2}|dx_1\ldots dx_k|^{\fr12}.$
The operator $W^{-1} TW$ will be an invariant pseudodifferential 
operator of the same order 
as $T,$ whose Schwartz kernel $K_{W^{-1}TW}$ has is given by
\begin{equation}
\label{propr2}
	K_{W^{-1}TW}(x,y,x' ,y' ) = A(e^{x'-x},y,y').
\end{equation}
and hence has support
\begin{equation}
\label{propr1}
\supp(K_{W^{-1}TW})\subset \R^k \times L \times \R^k \times L 
\end{equation}
where $L\subset \R^{n-k}$ is the compact set considered above. The distribution
$K_{W^{-1}TW}$ is a function ouside the diagonal satisfying the
estimates
\begin{equation}
\label{propr3}
	|K_{W^{-1}TW}(x,y,x' ,y' )| = | A(e^{x'-x},y,y')| 
   	\leq C_n e^{-n|x'-x|}\leq C_n (1 + |x'-x|)^n 
	\end{equation}
for $|x'-x|\geq 1$.

Consider the normalized Fourier transformation $\mathcal F$ in the first 
$k$-variables
$$
\gathered
	\mathcal F:\ll (\R^n) \lra \ll (\R^n) \\
	(\mathcal F v) (\xi, y) = (2\pi)^{-\fr{k}2} \int \,
	e^{-i \xi x} {v(x,y) dx}
\endgathered
$$

The action of $T_1 = \mathcal F^{-1} W^{-1} TW\mathcal F$ is given
explicitely by:
$$
	(T_1v) (\xi, y) = \int_{\R^{n-k}} \hat{T} 
	(\xi, y, y' ) v (\xi, y' ) dy' 
$$
where 
$\hat{T} (\xi, y, y' ) = \int \, e^{i \xi x} k_0 (x,y,y' ) dx$ 
(in distribution sense) and $k_0(x,y,y')=A(e^{x'-x},y,y')$. 
Standard estimates now show:
\begin{enumerate}
\item   $\hat{T}$ is differentiable of arbitrary order 
in $\xi$ because $k_0$ is quickly decreasing in $x$, equation \eqref{propr3}.

\item  If $T$ is of negative order then for each fixed $\xi$ the operator 
$\hat{T} (\xi)$ with Schwartz kernel $\hat{T}(\xi, y, y' )$ acting on $\ll (\R^{n-k})$
is a pseudodifferential operator of negative order with
compactly supported distribution kernel and hence it is a 
compact operator.

\item The function $D_{\xi_i}\hat{T}(\xi)$ is of the same form,
$D_{\xi_i}\hat{T}(\xi)=A_{T_i}(\xi)$ where $T_i=[x_i,T]$ is an operator
of lower order (see to \eqref{automorphism} and bellow).

\item   If $T$ is of order $< -n-1,$ then $\hat{T}$ is $C^1$ and
we have an estimate of the form
$$
	|\hat{T} (\xi,y,y' )| \le C (1+ |\xi|)^{-1}, \quad C>0.
$$
\end{enumerate}
Suppose $\operatorname{ord} (T) < -n-1.$ Since $\hat{T} (\xi, y, y')  = 0$ for $y$ or  
$y' $ outside a fixed compact set and $\hat{T}$ is of $C^1$-class, we 
obtain that $\|\hat{T}(\xi)\| \le C'  (1+ |\xi|)^{-1},$ for some uniform
constant $C' > 0$. This shows that the map
$$
q : \psi_m \Inal{-1}F\psi_m \lra 
C^{\infty}(\R^k, {\mathcal K})\,,\;q (T) = \hat{T}
$$
maps $\psi_m \Inal{-n-2}F\psi_m$ to functions 
vanishing in norm at infinity. We extend now the vanishing at
infinity property to all of $\psi_m \Inal{-1}F\psi_m$.

Let $T\in \psi_m \Inal{-1}F \psi_m.$ By construction 
$q((T^* T)^ l) (\xi) = (q(T)(\xi)^*q(T)(\xi))^l$ which, according to 
the previous arguments, is a smooth function with values compact 
operators and vanishing at infinity.  We then obtain
$$
\|q((T^* T)^ l) (\xi)\| = \|(q(T)(\xi)^*q(T)(\xi))^l \|
=\|(q(T)(\xi)^*q(T)(\xi)) \|^l=\|q(T)(\xi) \|^{2l}.
$$

If $2l > n +1 $ this shows  in view of the preceeding discussion
that  $\|q(T)(\xi)\|=\|\hat{T} (\xi)\|\to 0$ 
for $|\xi|\to \infty$. We thus obtained $q(\psi_m \Inal{-1}F \psi_m)
\subset C_0 (\R^k, \mathcal K)$. The density of $\cup_m \psi_m \Inal{-1}F \psi_m$
in ${\mathfrak K}_F$, lemma \ref{lemma3}, proves
$$
	q({\mathfrak K}_F) \subset C_0 (\R^k, \mathcal K).
$$

Let  $f \in C_0 (\R^k, \mathcal K)$, $f (\xi) = \hat{f}_0 (\xi) P_0$   
where $f_0 \in C^\infty_c (\R^k)$ is compactly supported,
$P_0 \in \Psi_c^{-\infty}(\R^{n-k})$ is a regularizing operator with
compactly supported Schwartz kernel and $\hat{f}_0$ is the
Fourier transform of $f_0$. Define $T \in \cup_m \psi_m \Inal{-\infty}F \psi_m$ 
to correspond to the kernel
$$
	A(x, y, x', y' ) = \hat f_0 ( - \ln x') P_0 (y, y' ).
$$
Then we have $q (T) = f.$ Since the functions $f$ as above form
a dense subset of $C_0 (\R^k, \mathcal K)$ the result follows.
\end{proof}

\begin{corollary} \label{SubDen}
The space $\cup_m \psi_m \Inal{-\infty}{F}\psi_m$ is dense in 
${\mathfrak K}_F$.
\end{corollary}

\begin{proof}The proof of the above proposition shows that 
$q(\cup_m \psi_m \Inal{-\infty}{F}\psi_m)$ is dense in 
$C_0(\RR^k, {\mathcal K}_F)$. Since $q$ is an isometric embedding the 
result follows from the above result.
\end{proof}

\begin{proposition} \label{SymbolMap2}
The symbol map \eqref{bExtPseudoSymbol}
gives a short exact sequence 
\begin{equation}
0\longrightarrow\bnInal{F}\longrightarrow \bInal{F}
\overset{\sigma_{0F}}
{-\!\!\!\!\longrightarrow}C(\bS^*M\vert_L)\longrightarrow 0
\end{equation}
where $\bnInal{F}$ is the norm closure of $\Inal{-1}F.$
\end{proposition}

\begin{proof} We know that the principal symbol is surjective,
Lemma \ref{SymbolMap}, so
we only need  to prove exactness in the middle. Observe that
$\bInal{F}^{(-1)} \subset \ker \ssigma_{0,F}$ by the continuity
of $\ssigma_{0,F}$. In order to prove the opposite inclusion
we need to show that if $T = \lim T_n,$ $T_n \in \Inal{0}F$ 
and $\lim \ssigma_{0,F}(T_n)=0$
then we can find a  sequence $T'_n$ with with the same limit as
$T_n$ and $\ssigma_{0,F}(T'_n)=0.$ We will prove this by induction on 
how singular the face $F$ is. If $F$ is a smooth manifold
without corners then the statement follows from 
lemma  \ref{lemma1} (it is a classical result).
If $F$ has `corners' we can assume by induction and using 
Corollary \ref{Cor1} that $\ssigma_{0,F}(T)=\ssigma_{0,F}(T_n)=0$ on any 
boundary hypersurface of $F.$ By the same argument
as in Lemma \ref{lemma3} (using the function $f$ that vanishes exactly 
on the union of all boundary hypersurfaces  of $F$ and the exhaustion
sequence $\psi_m$) we can further assume that all $T_n$ vanish 
in a neighborhood of the boundary.
This again reduces the proof to the case of a smooth manifold to which
the lemma \ref{lemma1} applies.
\end{proof}

We summarize our results in the following theorem.
Denote by $\mathcal K({\mathcal H})$  the algebra of compact operators 
on a Hilbert space ${\mathcal H}$, ${\mathcal K}_F={\mathcal K}(\ll (F))$
and by $\mathcal F_{l}(M)$
the set of faces of dimension $l$ of the manifold with corners $M.$

\begin{theorem}\label{Theorem.CompositionSeries} The norm closure
${\mathfrak A}_M$ of the algebra of order zero $b$-pseudodifferential 
operators  on the compact connected manifold with corners 
$M$ has a composition series consisting of closed ideals,
$$
	{\mathfrak A}_M \supset {\mathfrak I}_0 \supset {\mathfrak I}_1 
	\supset \ldots \supset {\mathfrak I}_n\,,\; n=\dim(M),
$$
where ${\mathfrak I}_0$ is the norm closure of $\overline {\Psi^{-1}_b}(M)$ and  
${\mathfrak I}_l$ is the norm closure of the ideal of order $-1,$ 
$b-$pseudodifferential operators on $M$ whose indicial parts vanish on all 
faces of dimension less than $l.$ The partial quotients are
determined by 
$\ssigma_{0}:{\mathfrak A}_M / {\mathfrak I}_0\overset{\sim}{\lra} C_0 (S^*M),$ 
and
$$
{\mathfrak I}_{l}/{\mathfrak I}_{l+1} \simeq
\bigoplus\limits_{F\in \mathcal F_{l} (M)}
C_0 (\R^{n-l}, {\mathcal K}_F)\, , \;\; 0\leq l\leq n .
$$
The composition series and the isomorphisms are natural with respect to 
maps of manifolds with corners which are local diffeomorphisms.
\end{theorem}

The last isomorphism reads  
${\mathfrak I}_n = {\mathcal K}_M$ for $l=n$ and 
$$
{\mathfrak I}_0/{\mathfrak I}_1 = 
\bigoplus\limits_{F\in \mathcal F_0 (M)} 
C_0 (N^*F)\simeq\bigoplus\limits_{F\in \mathcal F_0 (M)} C_0 (\R^n)
$$
for $l=0$ since ${\mathcal K}_F =\CC$ if $F$ has dimension $0$.

\begin{proof} The fact that the principal symbol induces an isomorphism 
$\sigma_0:{\mathfrak A}_M/{\mathfrak J}_0 \simeq C_0(\bS^*M)$
was proved in proposition \ref{SymbolMap2}.

Consider the morphism
$\Inop_{l}:\bInal{M} \to \oplus _{F\in \mathcal F_{l} (M)} \bInal{F}$.
The kernel of this morphism is by definition ${\mathfrak I}_{l+1}$. Moreover
$\Inop_l$ maps ${\mathfrak I}_l$ to 
$\oplus_F  {\mathfrak K}_F$, $F\in \mathcal F_{l} (M)$, by the definition
of ${\mathfrak K}_F$. This gives us an inclusion
$$	
        {\mathfrak J}_l/{\mathfrak J}_{l+1} \subset 
	\bigoplus\limits_{F\in \mathcal F_{l} (M)}{\mathfrak K}_F . 
$$ 
Corollary \ref{Cor1} tells us that this inclusion is actually an equality.
From proposition \ref{Prop5} we obtain the isomorphism
$\oplus_F{\mathfrak K}_F \simeq \oplus_{F}
C_0 (\R^{n-l}, {\mathcal K}_F)$, $F\in \mathcal F_{l} (M)$.

The naturality of the composition series follows from the naturality of
the principal symbol and of the indicial maps.
\end{proof}

The indicial algebras $\bInal{F}$ have similar composition series 
which are compatible with the indicial morphisms.

\begin{theorem}\label{Theorem.CompositionSeries2} The algebra
${\mathfrak A}_F$ has a composition series
$$
	{\mathfrak A}_MF\supset {\mathfrak J}_0 \supset {\mathfrak J}_1 
	\supset \ldots \supset {\mathfrak J}_n\,,\; n=\dim(F),
$$
where ${\mathfrak J}_0 =\overline {\Psi^{-1}_{b,I}}(N_+F)=\ker \sigma_{0F}$ and 
${\mathfrak J}_l$ is the closure of the ideal of order $-1,$ 
$b-$pseudodifferential operators whose indicial parts vanish on all 
faces $F' \subset F$ of dimension less than $l.$ The partial quotients are
determined by the natural isomorphisms 
$\ssigma_{0}:{\mathfrak A}_F/ {\mathfrak J}_0\overset{\sim}{\lra} 
C_0 (S^*M\vert_F),$ 
and
$$
{\mathfrak J}_{l}/{\mathfrak J}_{l+1} \simeq
\bigoplus\limits_{F'\in \mathcal F_{l} (F)}
C_0 (\R^{n-l}, {\mathcal K}_{F'})\, , \;\; 0\leq l\leq n .
$$
\end{theorem}

\begin{proof} The proof consists of a repetition of the arguments in the
proof of the preceeding theorem, replacing $M$ by $F$.
\end{proof}

\begin{corollary} \label{RefDense} 
The subalgebra $\Inal{-\infty}F$ is dense in $\bnInal{F}$.
\end{corollary}

\begin{proof} The inclusion $\Inal{-\infty}F \subset \bnInal{F}$
preserves the natural composition series of both algebras. By induction
it is enough to show that the subquotients of the first algebra are
dense in the second algebra. The subquotients of the second
algebra are (direct sums of) $\mathfrak K_F$'s, whereas the coresponding
subquotients of the first algebra contain the spaces 
$\cup_m \psi_m \Inal{-\infty}{F}\psi_m$ of corollary \ref{SubDen}. The 
result then follows from the same corollary.
\end{proof}

We have the following generalization of Proposition \ref{Prop5}.

\begin{corollary} \label{Cor.Suspension}
The  Fourier transform in the fiber direction 
establishes an isomorphism $\bnInal{F,M} \simeq C_0(\R^k,\bnInal{F,F})$
if the face $F$ has codimension $k$ in $M$.
\end{corollary}

\begin{proof} The map $q$ of Proposition \ref{Prop5} extends to a map
$$
	q:\Inal{(-n-1)}{F} \lra C^\infty(\R^k,\Inal{(-n-1)}{F_0}).
$$
This map is compatible with the composition series of
$\bInal{F}$ and $\bInal{F_0}$ of the above Theorem
and induces an isomorphism on the partial quotients (after completing in norm).
The density property in the above corollary completes the proof.
\end{proof}

\begin{corollary}\label{Cor.Prod}
If $M_1$ and $M_2$ are two manifolds with corners
then 
$$\bnInal{F_1 \times F_2,M_1 \times M_2} 
\simeq \bnInal{F_1,M_1} \otimes_{min} \bnInal{F_2,M_2}.$$
\end{corollary}

The tensor product $\otimes_{min}$ is the minimal tensor product of
two $C^*$-algebras and is defined as the completion in norm of
$\bnInal{M_1} \otimes \bnInal{M_2}$ acting on 
$\ll(M_1 \times M_2)$. (The space $\ll(M_1 \times M_2)$ is the
Hilbert space tensor product $ \ll(M_1)\hat\otimes\ll(M_2)$, it is
the completion  of the algebraic tensor product $ \ll(M_1)\otimes\ll(M_2)$ 
in the natural Hilber space norm.)

\begin{proof} We will assume that $F_i=M_i$, the general case being
proved similarly.
We have $\Psi_b^{-\infty}(M_1) \otimes \Psi_b^{-\infty}(M_2)
\subset \Psi_b^{-\infty}(M_1 \times M_2)$.
From  corollary \ref{RefDense} we conclude the existence of a morphism 
$\chi:\bnInal{M_1} \otimes_{min} \bnInal{M_2} \lra \bnInal{M_1 \times M_2}$
which preserves the composition series.  Moreover by direct inspection
the morphisms induced by $\chi$ on the subquotients are isomorphisms. 
If follows that $\chi$ is an isomorphism as well.
\end{proof}

For a compact manifold with boundary the results can be made even more
explicit. Using the above notation we have

\begin{theorem} If $M$ is a compact manifold with boundary then
$$
{\mathfrak I}_0 
= {\mathfrak I}_{n-1},\ {\mathfrak I}_{n-1}/{\mathfrak I}_n \simeq C_0
(\R,{\mathcal K}_{\pa M})$$
and  $ {\mathfrak A}_M/{\mathfrak I}_0 = C_0 (\bS^* M).$
The algebra $Q_M = {\mathfrak A}_M/{\mathfrak I}_n$ has the following
fibered product structure 
$Q_M \simeq Q_0 \subset C_0(\bS^*M) \oplus {\mathfrak A}_{\partial M},$ 
$Q_0 = \{ (f, T),\,f\vert_{\pa M} = \ssigma_{0,\partial M} (T)\}$.
The indicial algebra of the boundary, ${\mathfrak A}_{\partial M},$ 
fits into an exact sequence 
$$
	0 \lra C_0 (\R, {\mathcal K}) \lra {\mathfrak A}_{\partial M} 
	\overset{\sigma_{0\partial M}}{-\!\!\!-\!\!\!\longrightarrow} 
        C_0(\bS^*M|_{\partial M}) \to 0 \,.
$$
\end{theorem}

Let us take a closer look at the regularizing operators 
on $M_0=[0,\infty) \times \R^{n-1}.$ These operators have the form
\ref{SimpleAction} have the form
\begin{multline}
Tu(x,y_1,\ldots,y_{n-1})=\\
\int\limits_0^\infty\dots\int\limits_0^\infty
A(x,y_1,\ldots,y_{n-1},x',y'_1,\ldots,y'_{n-1})
u(x_1x'_1,y'_1,\ldots,y'_{n-1})
\frac{dx'}{x'}dy'_{1}\ldots dy'_{n-1}
\label{SimpleAction3}
\end{multline}
with $A$ a smooth compactly supported function in $(y,y')$ on 
$[0,\infty) \times \R^{n-1} \times [0,\infty) \times \R^{n-1}$, vanishing 
to all orders at $x'=0$. Using this description of operators we obtain
\begin{equation}
\label{eqWH} 
	\Psi^{-\infty}_b (M_0) \simeq  \Psi^{-\infty}_b (H_0)
   	\hat{\otimes} \mathcal R.
\end{equation}
where $H_0=[0,1\infty)$ and 
$\mathcal R = \Psi^{-\infty}_c (\R^{n-1})$ is the algebra of smoothing 
operators with compactly supported Schwartz kernels. The operators acting 
on $H=[0,1)$ can be described in a similar way using kernels, or can
be obtained as operators of the form $\phi T \phi$ where $\phi$
is a smooth, compactly supported function on  $[0,1)$ and $T$ is an
operator on $[0,\infty)$. Recall that we require all operators to 
have compactly supported distribution kernels.
The explicit form of the kernels in equation \eqref{SimpleAction3} 
for $n-1=0$ shows that $\Psi_b^{-\infty}(H)$ can be 
identified with a dense subalgebra of the algebra of Wiener-Hopf operators
on $[0,\infty)$. The algebra  ${\mathcal W}$ of Wiener-Hopf operators
on $[0,\infty)$ is defined to be the norm closed
algebra generated by the operators $T_\phi$ on $\ll([0,\infty))$ 
$$
(T_\phi) f(t)=\int_0^\infty \phi(t-s)f(s) ds\;,\;\; \phi \in \CIc(\RR).
$$
\begin{lemma} \label{lemmaWH} 
The substitution $x=e^{-t}$ identifies $\Psi_b^{-\infty}(H)$, $H=[0,1)$,
with a dense subalgebra of ${\mathcal W}$ containing the operators
$T_\phi$ defined above.  The indicial parts of these operators is given by
$\Inop_{\pa H,H}(T_\phi)=\phi \in {\mathfrak K}_{\pa H} = C_0(\RR)$.
\end{lemma}

\section{Computation of the $K$-groups}

Our starting point for the computation of the $K$-groups of the algebras
discussed in the previous section is the short exact sequence, of $C^*$-algebras,
$$
	0\lra {\mathcal K} \lra {\mathfrak A}_M \lra Q_M \lra 0.
$$
where $\mathcal K = \mathcal K(\ll(M))$. This exact sequence gives rise 
to the fundamental six-term exact sequence in $K$-theory (see \cite{Blackadar1})
\begin{gather}
\begin{CD}
	K_0(\mathcal K)  @>>> K_0 ({\mathfrak A}_M) @>>> K_0(Q_M) \\
	@A\partial AA	@.	@VV0V \\
	K_1(Q_M) @<<< K_1({\mathfrak A}_M) @<<< K_1(\mathcal K) 
\end{CD}\label{6term}
\end{gather}
where we are particularly interested in the $K$-groups of $Q_M.$ Here
$K_0(\mathcal K) \simeq \mathbb Z,$ and $K_1(\mathcal K) \simeq 0,$ so the
right vertical map is zero.

The left vertical arrow represents ``the index
map''. If $P$ is an $m\times m$ matrix with values 
b-pseudodifferential operators on $M$ which is fully elliptic, 
in the sense that  its image $u_P$ in $M_{m} (Q_M)$ is
invertible, and hence defines an element $[u_P]\in K_1(Q_M),$ then
\begin{equation}
        \label{Index}
	\partial [u_P] = \Ind(P) = \dim \ker P - \dim \ker P^* \in 
	\ZZ \simeq 
     K_0(\mathcal K)
\end{equation}
see \cite{Blackadar1,Connes1,Kasparov1}.
It will be useful for us to study the exact sequences
$$
	0\lra  {\mathfrak I}_{l}/{\mathfrak I}_{l+1} \lra 
	{\mathfrak I}_{l-1}/{\mathfrak I}_{l+1} \lra 
	{\mathfrak I}_{l-1}/{\mathfrak I}_{l} \to 0
$$
corresponding to the composition series described in Theorem
\ref{Theorem.CompositionSeries}.
We know that 
\begin{equation}
K_i (C_0 (\R^j, {\mathcal K})) \simeq \cases 
\mathbb Z \quad \text{if } i+j \text{ is even} \\ 0 \quad \text{otherwise.}
\endcases\label{KPeriodicity}
\end{equation}
We shall fix these isomorphisms uniquely as follows.  For $j=0,$ $K_0
({\mathcal K}) \overset{\sim}{\lra}\mathbb Z$ will be the dimension
function; it is induced by the trace. For $j>0$ we define the
isomorphisms in \eqref{KPeriodicity} by induction to be compatible
with the isomorphisms
$$
\gathered
	\mathbb Z \simeq K_{2l-j+1} (C_0 (\R^{j-1}, \mathcal K)) 
	\overset{\partial}{\lra}  \\
	\lra K_{2l-j} (C_0((0,\infty)\times \R^{j-1}, \mathcal K)) 
	\overset{\sim}{\lra} K_{2l-j} (C_0(\R^j, \mathcal K))
\endgathered
$$
where the boundary map corresponds to the exact sequence of 
$C^*$-algebras
$$
	0\lra C_0 ((0, \infty) \times \R^{j-1}, \mathcal K) \lra 
	C_0 ([0, \infty) \times \R^{j-1}, \mathcal K) \lra 
	C_0 (\R^{j-1}, \mathcal K) \to 0.
$$
Denote, for any $C^*$-algebra $A$ by $SA = C_0 (\R, A) = C_0 (\R) 
\otimes_{\min} A,$ $S^k A = C_0 (\R^k, A).$
Define  
$F'_0 = \{(0\cdc0,0)\} \times \R^{l-1}, F_0 = \{(0\cdc 0)\} \times 
[0,1) \times \R^{l-1}$ and $M_0=[0,1)^{n-l+1} \times \R^{l-1}$, 
$F'_0 \subset F_0\subset M_0.$  Also let $H= [0,1)$ and 
$0\subset {\mathcal L}_1 \subset {\mathcal L}_0 \subset {\mathfrak A}_H$ be the 
canonical composition series of ${\mathfrak A}_H,$ such that ${\mathcal L}_0 = 
{\mathfrak A}^{(-1)}_H,$ ${\mathcal L}_0/{\mathcal L}_1 \simeq C_0 (\R),$ 
${\mathcal L}_1 \simeq \mathcal K$, see Theorem \ref{Theorem.CompositionSeries}.

\begin{lemma} Denote ${\mathcal K}_1=\mathcal K(\ll(\RR^{l-1}))$. 
There exists a commutative diagram 
$$
\begin{CD}
	0 @>>> \ker (\Inop_{F'_0,F_0}) @>>> {\mathfrak A}^{(-1)}_{F_0,M} @>>>
	 {\mathfrak A}^{(-1)}_{F'_0,M} @>>> 0 \\
	@.	@VVV	 @VVV		@VVV \\
	0 @>>> S^{n-l} {\mathcal L}_1 \otimes {\mathcal K}_1 
        @>>> S^{n-l} {\mathcal L}_0 \otimes 
	\mathcal K_1 @>>> S^{n-l+1} \mathcal K_1  @>>> 0
\end{CD}
$$
in which all vertical arrows are isomorphisms, the bottom exact 
sequence is obtained from $0 \to {\mathcal L}_1 \to {\mathcal L}_0 
\to C_0 (\R) \to 0$ by 
tensoring with $C_0 (\R^{n-l}, \mathcal K_1)$ and the boundary map
$$
\gathered
	\partial: K_{n-l+1} (C_0 (\R^{n-l+1}, \mathcal K_1)) \simeq 
	K_{n-l+1} ({\mathfrak A}_{F'_0}^{(-1)}) \lra \\
 	K_{n-l} (\ker (\Inop_{F'_0,F_0}))
	\simeq K_{n-l} (C_0 (\R^{n-l}, \mathcal K_1))
\endgathered
$$
is (the inverse of) the canonical isomorphism.
\end{lemma}

\begin{proof} Let $H^{l-1}=[0,1)^{l-1}$ and $F_1=\{(0,\ldots,0)\} \times H$,
$F_1 \subset H^{l-1}$.
It follows from the corollary \ref{Cor.Prod} that the algebra 
${\mathfrak A}^{(-1)}_{F_0,M_0}$ is isomorphic to 
$\bnInal{F_1,H^{l-1}} \otimes \bnInal{\RR^{l-1}}$.
Moreover $\bnInal{\RR^{l-1}} ={\mathcal K}(\ll(\RR^{l-1}))$.
The corollary \ref{Cor.Suspension} further gives
$\bnInal{F_1,H^{l-1}} \simeq C_0(\RR^{n-l},\bnInal{F_1,F_1})
=C_0(\RR^{n-l},\bnInal{H})$. The first commutative diagram then is just
an expression of the composition series of $\bnInal{H}$,
Theorem \ref{Theorem.CompositionSeries}. The lemma \ref{lemmaWH}
reduces the computation of the connecting morphism $\pa$ to that of the
connecting morphism of the Wiener-Hopf exact sequence
(\ie the Wiener-Hopf extension). This is a well known and easy computation. It 
amounts to the fact that the multiplication by $z$ has index $-1$ on
the Hardy space $H^2(S^1)$ of the unit circle $S^1$. See \cite{Blackadar1}
for more details.\end{proof}

From Theorem \ref{Theorem.CompositionSeries} we then know that 
$$
	K_i ({\mathfrak I}_{l}/{\mathfrak I}_{l+1}) \simeq 
	\cases\bigoplus\limits_{F\in \mathcal F_{l} (M)} \mathbb Z 
	\quad \text{if } n-l+i \quad \text{is even} \\
	0 \qquad\qquad \text{otherwise}
	\endcases
$$
where $n=\dim M.$
Fix from now on an orientation of each face of $M,$ including 
$M$ itself.  No compatibilities are required.  This will uniquely 
determine the above isomorphisms.

Define the incidence number $[F: F']$ for any codimension 1 face 
$F'$ of $F$ as usual.  It depends on the choice of 
orientations of $F$ and $F'.$  Given an orientation $e_1\cdc e_k \in 
T_p F,$ consisting of linearly independent vectors, 
we can assume that $p \in F',$ $e_1$ 
points inward and $e_2\cdc e_k \in T_p F'.$  Then we obtain a new orientation of $F'$ 
using $e_2\cdc e_k.$  The incidence number $[F: F']$ will be $+1$ 
if this new orientation of $F'$ coincides with the old one, $-1$ 
otherwise. Define $[F:F']=0$ if $F'$ is not a boundary hypersurface of $F$.

\begin{theorem} Suppose $n-l+i$ is even. Then the matrix of the boundary map
$$
	\partial: K_{i-1} ({\mathfrak I}_{l-1}/{\mathfrak I}_{l}) \simeq 
	\bigoplus\limits_{F' \in \mathcal F_{l-1} (M)} \mathbb Z 
	\lra\bigoplus\limits_{F\in \mathcal F_{l} (M)} \mathbb Z 
	\simeq K_i ({\mathfrak I}_{l}/{\mathfrak I}_{l+1}) 
$$
is given by the incidence matrix. If $n-l +i$ is odd, then $\partial=0.$
\end{theorem}

\begin{proof} Denote by $e_{F'} \in K_{i-1} ({\mathfrak
I}_{l-1}/{\mathfrak I}_l),$ $e_F \in K_i ({\mathfrak
I}_l/{\mathfrak I}_{l+1})$ the canonical generators of these groups.
We need to show that
$$
	\partial (e_{F'}) = \sum\limits_{F\in \mathcal F_l(M)} \,
	[F:F'] e_F.
$$
The idea of the proof is to reduce the computation to the case 
$M=M_0$, $F=F_0$ and $F'=F'_0$ considered in the preceeding lemma
$$
M_0= H^{n-l+1} \times \R^{l-1}\, , \;\; F_0 = 
\{(0\cdc 0)\} \times H \times   \R^{l-1}\, , \;\; F'_0 = 
\{(0\cdc 0,0)\} \times   \R^{l-1}
$$ 
and $F'=F'_0$ the face of minimal dimension.

Choose a point  $p\in F'.$ There exists a diffeomorphism 
$\varphi: M_0 \lra M$ of manifolds with corners onto an open 
neighborhood of $p$ such that $p \in \varphi(F_0')$.  
Since we considered only pseudodifferential 
operators with compactly supported Schwartz kernels, the open map 
$\varphi$ induces an inclusion
$$
	\varphi_0: {\mathfrak A}_{M_0} \lra {\mathfrak A}_M  
     \qquad(\text{label{ eq. parts}})
$$
which commutes with the indicial maps and hence preserves the 
canonical composition series of the Theorem \ref{Theorem.CompositionSeries}:  
$
\varphi_0 ({\mathfrak I}^{(0)}_{l}) \subset {\mathfrak I}_{l}
$ where 
${\mathfrak I}^{(0)}_{n} \subset {\mathfrak I}^{(0)}_{n-1} \subset \dots \subset 
{\mathfrak A}_{M_0}$ is the composition series associated to $\bInal{M_0}.$  
This composition series has the following properties
$$
	{\mathfrak I}^{(0)}_{l-1} =\dots = {\mathfrak I}^{(0)}_0=
\bnInal{M_0}, \quad 
	{\mathfrak I}^{(0)}_{l-1}/ {\mathfrak I}^{(0)}_l \simeq C_0 
	\big(\R^{l-1}, {\mathcal K}(\ll(\RR^{l-1})\big)
$$
and the induced map on $K$-theory 
$$
	\mathbb Z e_{F'_0} = K_{i-1} ({\mathfrak I}^{(0)}_{l-1}/ 
	{\mathfrak I}^{(0)}_l) \lra K_{i-1} 
	({\mathfrak I}_{l-1}/ {\mathfrak I}_l)
$$
maps  $e_{F'_0}$ to $e_{F'}.$  We will first compute the boundary map
$$
	\partial_0:   K_{i-1} ({\mathfrak I}^{(0)}_{l-1}/ 
	{\mathfrak I}^{(0)}_l) \lra K_i ({\mathfrak I}^{(0)}_l/ 
	{\mathfrak I}^{(0)}_{l+1})= \oplus 
	\mathbb Z e_{F^{(0)}_j}.
$$
with the generators of the second group being indexed by the faces 
$F^{(0)}_j$ of $M_0$ 
dimension $l.$ Since the boundary map in $K$-theory is natural, 
$\varphi_* \pa_0=\pa \varphi_*$, 
and $\varphi_{0*} (e_{F^{(0)}_j}) = e_{F_j}$ if 
$\varphi (F^{(0)}_j) \subset F_j,$ the boundary morphism $\partial_0$
will determine $\partial$ and thus prove the theorem.  Label the faces 
$F^{(0)}_0\cdc F^{(0)}_{n-l}$ in the order given by the additional 
coordinate (thus $F_0^{(0)}=F_0$). 
It is enough to compute $\partial_0$ for an arbitrary 
choice of orientations, so we can choose the canonical one (given by 
the order of components).  We then need to prove that the 
coefficients $c_j$ defined by 
$\partial e_{F'_0} = \sum c_j e_{F^{(0)}_j}$ 
satisfy $c_j = (-1)^{j}.$  By symmetry it is enough to assume $j=0.$  
The indicial map $\Inop_{F_0,M_0}$  restricts to  an onto morphism 
$$
	\psi=\Inop_{F_0,M_0}:\bnInal{M_0}=
	{\mathfrak I}^{(0)}_{l-1} \lra {\mathfrak A}^{(-1)}_{F_0} = 
	\ker \ssigma_{0,F_0} \;\;\;\; 
$$
such that 
$
	\psi({\mathfrak I}^{(0)}_l) = \ker (\Inop_{F'_0,F_0} : 
	{\mathfrak A}^{(-1)}_{F_0} \to {\mathfrak A}^{(-1)}_{F'_0}) 
	\simeq C_0 (\R^{n-l+1}, \mathcal K)
$ and $\psi({\mathfrak I}_{l +1}^{(0)})=0.$ (Recall that
$F_0 = F^{(0)}_0, F'_0 \subset F_0$.) The induced morphism
$$
	\psi_*: K_i ({\mathfrak I}^{(0)}_l / {\mathfrak I}^{(0)}_{l+1}) \to \mathbb Z
$$
is the projection onto the first component (\ie it gives the coefficient
of $e_{F_0}$).  Using again the 
naturality of the exact sequence in $K$-theory we further reduce
the proof to the computation of the boundary map in the exact sequence
$$
	0 \lra \ker (\Inop_{F'_0,F_0})  \lra 
	\bnInal{F_0} \lra \bnInal{F'_0} \to 0
$$
This computation is the content of previous lemma.  This completes the
proof.\end{proof}

The above theorem identifies the $E^1$ term and the $d_1$-differential of
the spectral sequence associated to the above composition series by
C. Schochet \cite{Schochet1}.

The results we obtained on the structure of the norm closure of the 
algebras of b-pseudodifferential operators on a manifold with corners
extend immediately to fibrations an families of pseudodifferential
operators acting on the fibers.
We state these results for  families of manifolds with boundary. 

Let $\pi: Z\to X$ be a fiber bundle such that $X$ is a Hausdorf 
locally compact space and the fibers are $C^{\infty}$-manifolds with 
boundary with the smooth structure varying continuously in $x\in X.$

We consider the algebra $\Psi^0_{b,comp} (Z \to X)$ of families of 
b-pseudodifferential operators of order $0$ on the fibers of $Z \to X$ 
whose Schwartz kernels are globally 
compactly supported.  Denote its norm closure by ${\mathfrak A}_Z.$
The indicial operator 
$$
	\Inop_{\partial Z, Z}: \Psi^0_{b,comp} (Z \to X) \lra \Psi^0 
	(N_+ \partial Z \to X)^{\R^*_+}
$$
is defined analogously and extends to ${\mathfrak A}_Z.$  Denote by 
$\bS^*_f Z$ the $b$-cosphere bundle along the fibers of $Z \to X,$ by 
${\mathfrak I}_{n-1} = \ker \sigma_0 : {\mathfrak A}_Z \lra C_0 (\bS^*_f Z)$ and by 
${\mathfrak I}_n = \ker (\Inop_{\partial Z, Z}) \cap {\mathfrak I}_{n-1}.$

\begin{theorem}\label{Theorem.IndicialCompletions} For a fibration 
$Z\lra X$ with fibers manifolds with boundary and with ${\mathcal L}_0=\bnInal{H}$,
$H=[0,1)$ as above, we have  isomorphisms
${\mathfrak I}_n \simeq C_0 (X, \mathcal K),$
${\mathfrak I}_{n-1} =\overline\Psi_b^{-1}(Z\to X)\simeq C_0 (X, {\mathcal L}_0),$
${\mathfrak I}_{n-1}/{\mathfrak I}_n \simeq C_0 (X \times \R, \mathcal K)$
and
${\mathfrak A}_Z/ {\mathfrak I}_{n-1} \simeq C_0 (\bS_f^* Z).$
\end{theorem}

\begin{proof} Just repeat all the arguments above using an extra parameter 
$x\in X.$\end{proof}

The above theorem allows us to determine very explicitly the $K$-theory 
groups of the norm-closed algebras associated to families of manifolds 
with boundary.

\begin{theorem} Let $Z \to X$ be as above and set $Q_Z = {\mathfrak
A}_Z/{\mathfrak I}_n$. Then the principal symbol $\sigma_0$ induces isomorphisms
\begin{equation}
	\ssigma_{0*}: K_i ({\mathfrak A}_Z) \simeq K_i (C_0 (\bS_f^* Z)) \simeq K^i 
	(\bS_f^* Z)
\label{symbolKisomorphisms}
\end{equation}
and the boundary map
$
	\partial: K_i ({\mathfrak A}_Z/ {\mathfrak I}_{n-1}) \lra K_{i+1} 
	({\mathfrak I}_{n-1}/ {\mathfrak I}_n)
$
is zero so there is a natural short exact sequence
\begin{equation}
	0 \lra K^{i} (\R \times X) \overset{j_Z}{\lra} K_i (Q_Z) \lra K^i 
	(\bS_f^* Z) \to 0
\label{jshortexact}
\end{equation}
\end{theorem}

\begin{proof} The groups $K_i (C (X, {\mathcal L}_0))$ can be computed
using the K\"unneth formula \cite{Blackadar1,Rosenberg-Schochet1}.  Since
$K_* ({\mathcal L}_0) = 0$ it follows that $K_* ({\mathfrak I}_{n-1})=0,$
which, in view of Theorem \ref{Theorem.IndicialCompletions}, proves the
first part of the theorem.

In order to prove \eqref{jshortexact} observe that, using
\eqref{symbolKisomorphisms}, the composite map
$$
	K_i ({\mathfrak A}_Z) \lra K_i (Q_Z)\lra K^i (\bS_f^* Z)
$$
is surjective, so the map from $K_i(Q_Z)$ to $K^i(\bS_f^* Z)$
in \eqref{jshortexact} is also surjective. This shows that $\partial = 0.$
\end{proof}

One important problem is explicitely compute the family index map 
\cite{Atiyah-Singer3} 
$$
	\Ind = \partial: K_i (Q_Z) \lra 
	K_{i-1}(C_0(X,{\mathcal K})) \simeq K^{i-1} (X)
$$
corresponding to the exact sequence
$
	0 \lra C_0 (X, {\mathcal K}) \lra {\mathfrak A}_Z \lra Q_Z \lra 0
$
along the lines of \cite{MelroseNistor2}. A consequence of our
computations is the following corollary.

\begin{corollary} The composition 
$$
	\Ind \circ j_Z:K^{i} (\R \times X) \lra K^{i-1} (X)
$$
is the canonical isomorphism.
\end{corollary}

\begin{proof} Indeed, $\Ind \circ j_Z$ is by naturality the
boundary map in the six term $K$-Theory exact sequence  associated
to the exact sequence 
$$
0 \lra {\mathfrak I}_{n-1} \lra {\mathfrak I}_{n} \lra 
{\mathfrak I}_{n-1}/{\mathfrak I}_{n} \lra 0.
$$
Since $K_i({mathfrak I}_n)=0$ it follows that the connecting (\ie boundary)
morphism in the above six term exact sequence is an isomorphism.  
The descriptions of the ideals ${\mathfrak I}_{n-1}$ and 
${\mathfrak I}_{n}$ in the previous theorem then completes the proof.
\end{proof}

\section{An $\R^q$-equivariant index theorem}

Consider a smooth manifold $X$ endowed with a proper, free action of
$\R^q.$ We can assume that $X=\R^q \times F$ with $\R^q$ acting by
translations. If $D$ is an elliptic differential operator on $\R^q\times F$
which is invariant under the action of $\R^q$ it is natural to look for
index type invariants of $D.$ We construct and compute such invariants using 
results from the previous sections. More generally we consider elliptic
matrices.

Denote by $n$ the dimension $F.$ We have already defined the algebra,
$\Inal0{F},$ of $\R^q$-invariant, pseudodifferential operators order
$0$ on $\R^q \times F$ when we studied the indicial algebra at a boundary face,
$F,$ of codimension $q$ of the noncompact manifold $M.$

Appealing to the same philosophy as before, we shall to consider
the closure in norm of this algebra, denoted $\bInal{F}$ as before.
We then know from Proposition \ref{SymbolMap2} that there is an exact sequence
$$
	0 \lra C_0 (\R^q, \mathcal K) \lra {\mathfrak A}_{F}
	\overset{\ssigma_{0F}}{-\!\!\!\lra} C_0 (S^* M\vert_{F}) \to 0.
$$
The connecting morphism (boundary map)
$$
	\partial:K_{q+1}(C_0 (\bS^* M\vert_{F})) = K^0(S^* M\vert_{F}) 
	\to K_q(C_0 (\R^q,\mathcal K))\simeq \mathbb Z
$$
can be interpreted as an $\RR^q$-equivariant index, and its computation 
will then be an $\R^q$-equivariant index theorem. Above we have 
used the standard isomorphism $K_q(C_0 (\R^q, \mathcal K))\simeq \mathbb Z$
as explained in section 2. 

Denote by $Y=S^*M\vert_{F}=S^*M/\RR^q$ and orient it
as the boundary of (the dual of) $\R^q \times TM$ if $TM$ is oriented as
an almost complex manifold (as in the Atiyah-Singer index theorem). Also
denote by ${\mathcal T}(F) \in H^{even}(F)$ 
the Todd class of the complexified cotangent bundle
of $F$ and by $Ch$ the Chern character.

\begin{theorem}
Let $a$ be an element of $K_{q+1}(C_0 (Y)),$ $Y=S^*M/\RR^q$.
Then
$$
	\partial (a)=(-1)^n \langle Ch(a)p^*{\mathcal T}(F), [Y] \rangle
$$
where $[Y]$ is the fundamental class of $Y$ oriented
as above, $p:Y\to F$ is the projection and $n=\dim F$.
\end{theorem}

\begin{proof}
We will prove the above theorem by induction on $q.$ Consider the manifolds with
corners
$$
	M=[0,1)^q \times F\; ,\;\; F_0=[0,1) \times \{(0\cdc 0)\} \times F \subset M
$$ 
and identify $F$ with
$\{(0,0\cdc 0)\} \times F \subset F_0.$ 
Let ${\mathfrak A}_{F_0}$ be as above and consider the following 
two onto morphisms and their kernels:
$$
\aligned
	\ssigma_{0,F_0}:{\mathfrak A}_{F_0} \lra C_0(Y) \;,\;\; I=\ker 
	(\ssigma_{0,F_0})
\\
	\Inop_{F,F_0}:{\mathfrak A}_{F_0} \lra  {\mathfrak A}_{F} \;,\;\; J=\ker 
	(\Inop_{F,F_0})
\endaligned
$$
Denote by $\partial_1=\partial,\partial_2,\partial_3,\partial_4$ the following
boundary maps in $K$-theory:
$$
\gathered
	\partial_1:K_{q+1}({\mathfrak A}_{F_0}/(I+J))=
	K_{q+1}(C_0 (Y)) \to K_q(I/I\cap J)
	=K_q(C_0 (\R^q,\mathcal K))\simeq \mathbb Z
\\
	\partial_2:K_{q+1}({\mathfrak A}_{F_0}/(I+J))
	\to K_q(J/I\cap J)=K_q(C_0(Y_0))=K^q(Y_0)
\\
	\partial_3:K_q(I/I\cap J) \to K_{q-1}(I\cap J)
	=K_{q-1}(C_0(\R^{q-1},\mathcal K))\simeq \mathbb Z
\\
	\partial_4:K_q(J/I\cap J) \to K_{q-1}(I\cap J)
	\simeq \mathbb Z
\endgathered
$$
where $Y_0=S^* M\vert_{F_0} \setminus S^* M\vert_{F}=(0,1)\times Y$ 
and we have used the determination of the partial quotients given by theorem
\ref{Theorem.CompositionSeries}.

We have that $\partial_3 \partial +\partial_4 \partial_2=0$ due to the fact that
the composition of connecting morphisms
$$
	K_{q+1}({\mathfrak A}_{F_0}/(I+J)) \longrightarrow
	K_q((I+J)/I\cap J) \longrightarrow K_{q-1}(I\cap J)
$$
is $0$ and $(I+J)/I\cap J \simeq I/I\cap J \oplus J/I\cap J .$ 

Since $\partial_3$ is the canonical isomorphism
we obtain after identifications $\partial=-\partial_4\partial_2$.
We have $\dim F_0=\dim F+1=n+1,$ hence the codimension $q$
decreases by $1$. By induction, the theorem is true for $F_0$:
$$
	\partial_4 (a)=(-1)^{n+1} \langle Ch(a)p_0^{*}{\mathcal T}(F_0), [Y_0] \rangle
$$
for all $a \in K_q(C_0(Y_0)),$. We have denoted by  $p_0:Y_0 \to (0,1) \times F$
the projection. Let $\pi:Y_0=(0,1) \times Y\to Y$ be the projection
onto the second component. Then $p_0^{*}{\mathcal T}(F_0)=\pi^*p^*{\mathcal T}(F).$
Denote by $\partial_c$ the boundary map
$H_c^*(Y) \to H_c^{*+1}(Y_0).$
Using the fact that the Chern character is compatible with
the boundary maps in $K$-theory (this fact is proved for algebras in general 
in \cite{Nistor2,Nistor1}) we obtain
$$
\aligned
	\partial (a) & =-\partial_4 \partial_2 (a)
	= -(-1)^{n+1} \langle Ch(\partial_2 (a))
	p_0^{*}{\mathcal T}(F_0), [Y_0] \rangle
\\
	&= (-1)^{n} \langle \partial_c Ch(a)\pi^*p^*{\mathcal T}(F), [Y_0] \rangle
\\
	&= (-1)^{n} \langle \partial_c (Ch(a)p^*{\mathcal T}(F)), [Y_0] \rangle
\\
	&= (-1)^{n} \langle Ch(a)p^*{\mathcal T}(F), [Y] \rangle
\endaligned
$$
where the last eqality is Stokes' theorem. The theorem is proved.
\end{proof}	

\begin{lemma}The connecting morphism $\pa$ in the previous theorem
is onto.
\end{lemma}

\begin{proof} Choose a small contractible open subset $U\subset Y$ and choose 
the class  $a  \in K_{q+1}(C_0 (Y))$ to come from a generator
of $K_{q+1}(C_0 (Y\vert_U))$. Then $\pa a = \pm 1$.
\end{proof}

The following corollary determines the $K$-theory groups of the 
`higher indicial algebras'. The space $\bS^*M\vert_F$ in the statement
of the corollary plays the role of $Y$ in the last theorem.

\begin{corollary} If $F$ is a smooth face of the manifold with
corners $M$, of codimension $q$, then the $K$-theory groups of $\bInal{F}$ are
$K_{q+1}(\bInal{F}) \simeq \ker(\pa : K_{q+1}(\bS^*M\vert_F) \to \ZZ)$
and $K_{q}(\bInal{F}) \simeq K_q(\bS^*M\vert_F)$.
\end{corollary}

\begin{proof} The result follows from the $K$-theory six term exact
sequence applied to the short exact sequence in \ref{SymbolMap2} and the 
determination of the connecting morphism of that exact sequence obtained 
in the above lemma.
\end{proof}

One should compare the above theorem with other equivariant index theorems
for {\it noncompact} groups \cite{Atiyah2,Singer2} for discrete groups and
\cite{Connes-Moscovici1} for connected Lie groups.

\section{Final comments}

Consider the case $q=1$ in the preceeding theorem, and let as above
$F$ be a smooth manifold (without corners). The 
indicial algebra ${\mathfrak A}_{F}$ fits into an exact sequence 
$$
	0 \lra C_0 (\R, \mathcal K) \lra {\mathfrak A}_{F}
	\overset{\ssigma_{0,F}}{-\!\!\!\lra} C_0 (S^* M\vert_{F}) \to 0.
$$
The above results imply that the morphism
$$
	K_1 (C_0 (\R, \mathcal K)) \lra K_1 ({\mathfrak A}_{F})
$$
vanishes.  However at the level of the {\em algebraic} $K_1$ and 
{\em uncompleted} algebras the morphism
$$
	i: K^{\text{alg}}_1 (\Psi^{-\infty}_b (\R \times F)^{\R}) \lra 
	K^{\text{alg}}_1 (\Psi^{\infty}_b (\R \times F)^{\R})
$$
is not zero.  This follows from results of Melrose \cite{Melrose46} 
who proves the existence of a onto morphism $\eta: K^{\text{alg}}_1 
(\Psi^{\infty}_b (\R \times F)^{\R}) \lra \mathbb C$ wich coincides with the usual 
$\eta$-invariant of \cite{Atiyah-Patodi-Singer1} for admissiable Dirac operators.
Moreover  the composition $\eta\circ i$ computes the spectral flow.

In  subsequent papers \cite{MelroseNistor2,MelroseNistor3}
we will use the above observation to exted the
result of \cite{Melrose46} to familes and to study the relation
between the $\eta$-invariant and cyclic cohomology \cite{Connes1}.

\begin{appendix}
\section*[Appendix]{Appendix. The algebra of b-pseudodifferential operators}

For completeness we give a direction definition of the algebra of
b-pseudodiff\-er\-ent\-ial operators from which all the basic properties
can be deduced.  As model space consider the product of intervals
$M=M_k=[-1,1]^k.$ The space $\Psi_b^m(M)$ may then be defined directly as a
space of kernels. They are also defined on a similar model space $M^2_{\bar
b}=[-1,1]^{2k},$ although this should {\em not} be thought of as the
product of $M$ with itself. Rather the smooth map
\begin{equation}
\gathered
\bar\beta:M^2_{\bar b}\supset(-1,1)^k\times[-1,1]^k\ni(\tau,R)\longmapsto
(x,x')\in M^2\\
x_i=\frac{2R_i+\tau+R_i\tau_i}{2+\tau_i+R_i\tau_i},\
x_i'=\frac{2R_i-\tau_i-R_i\tau_i}{2-\tau_i-R_i\tau_i}\\
\Longleftrightarrow 
R_i=\frac{x_i+x_i'-2x_ix_i'}{2-x_i-x_i'+x_ix_i'},\
\tau_i=\frac{x_i-x_i'}{1-x_ix_i'}
\endgathered\label{bIndentification}
\end{equation}
is used to identify the interiors of $M^2$ and $M^2_{\bar b}.$ Consider the
`diagonal' submanifold $\Diag_{\bar b}=\{\tau_i=0,i=1,\dots,k\}$ and the
boundary hypersurfaces $B_i^{\pm}=\{R_i=\pm1\}.$ Then as
a linear space
\begin{equation}
\Psi_b^m(M)=\big\{A=A'\nu,\ A\in
I^m(M_{\bar b}^2,\Diag_{\bar b});A\equiv0\text{ at }B_i^{\pm}\ \forall\
i\big\}.\label{bpseudosDefined}
\end{equation}
Here $\nu$ is a `right density' namely
$$
\nu=\frac{|dx'_1\dots dx'_k|}{(1-(x'_1)^2)\dots(1-(x'_k)^2)}
$$
and $I^m(M^2_b,\Diag_{\bar b})$ is the space of conormal distributions. If
the kernel space is embedded,
$M^2_{\bar b}=[-1,1]^{2k}\hookrightarrow \RR^{2k},$ and then rotated so that
the linear extension of $\Diag_{\bar b}$ becomes the usual diagonal then
this space is precisely the restriction to the image of $M^2_{\bar b}$ of
the space of kernels of (polyhomogeneous) pseudodifferential operators on
$\RR^k.$ Since these kernels are smooth away from $\Diag_{\bar b}$ the
condition in \eqref{bpseudosDefined} that the kernels vanish in Taylor
series at the boundary faces $B_i^{\pm}$ is meaningful.

The identification $\bar\beta$ in \eqref{bIndentification}
transforms the space $\Psi_b^m(M)$ to a subspace of the space of
(extendible) distributional right densities on $M^2,$ so by the Schwartz
kernel theorem each element defines an operator. If $u\in\dCI(M),$ the
space of smooth functions vanishing at the boundary in Taylor series, then
in principle $Au,$ for $A\in\Psi_b^m(M),$ is a distribution on the interior
of $M.$ In fact
$$
A:\dCI(M)\longrightarrow \dCI(M)\ \forall\ A\in\Psi_b^m(M).
$$
The space $\Psi_b^m(M)$ is a $\CI(M^2)$ module, where the smooth functions
on $M^2$ are lifted under $\bar\beta$ to (generally non-smooth) functions on
$M^2_{\bar b},$ it is also invariant under conjugation by a
diffeomorphism. As a consequence for a general manifold with corners the 
space $\Psi^m_b(X)$ can be defined by localization.

\begin{definition}\label{BasicDefinition} If $X$ is any manifold with
corners then the space $\Psi^m_b(X)$ consists of those operators
$A:\dCIc(X)\longrightarrow \dCI(X)$ such that if $\phi\in\CIc(X)$ has
support in a coordinate patch diffeomorphic to a relatively open subset of
$M_k$ then the image of the localized operator $\phi A \phi$ under the
diffeomorphism is an element of $\Psi_b^m(M_k)$ and if
$\phi,\phi'\in\CIc(X)$ have disjoint supports in (possibly different) such
coordinate patches then the image of $\phi' A\phi$ is in $\Psi^0_b(M_k).$
\end{definition}

The
identification of $[-1,1]$ as the projective compactification of $[0,\infty)$
\begin{equation}
[0,\infty)\ni t\longmapsto x=\frac{t-1}{t+1}\in[-1,1]
\label{ProjectiveCompactification}
\end{equation}
interprets $M$ as a compactification of $[0,\infty)^k.$ This induces an
action of $(0,\infty)^k$ on $M$ and the invariant elements of $\Psi_b^m(M)$
correspond exactly to the kernels which are independent of the variables $R_i.$
The compactification similarly reduces $M^2_{\bar b}$ to a
compactification of $[0,\infty)^{2k}.$ The action of $A\in\Psi_b^m(M)$ can
then be written
\begin{equation}
Au(t_1,\dots,t_k)=\int\limits_0^\infty\dots\int\limits_0^\infty
A(t_1,\dots,t_k,s_1,\dots,s_k)u(s_1t_1,\dots,s_kt_k)
\frac{ds_1}{s_1}\dots\frac{ds_k}{s_k}
\label{SimpleAction}
\end{equation}
where $A(t,s)$ is smooth in $t,$ conormal in $s$ at $s=1$ and vanishes
rapidly (uniformly on compact sets) with all derivatives as any $s_i\to0$
or $\infty.$

An alternative definition of $b$-pseudodifferential operators is
obtained by considering as in \cite{Melrose25} operators of the
form
$$
	(T u) (x, y) = (2\pi)^{-n} \int \,  
	e^{i (x-x') \xi + (y-y')\eta}a(x,y,x\xi,\eta)
	u(x' ,y' ) dx'  dy'  d\xi d\eta
$$
where $a(x,y,\xi,\eta)$ is a classical (i.e. 1-step polyhomogeneous 
of integral order) symbol satisfying a certain lacunary condition.
The operator $T$ is seen to be invariant if and only if it is of the form
$$
	(T u) (x, y) = (2\pi)^{-n} \int \,  
	e^{i (x-x') \xi + (y-y')\eta}  
	a(0,y,x\xi,\eta)
	u(x' ,y' ) dx'  dy'  d\xi d\eta.
$$
\end{appendix}


\providecommand{\bysame}{\leavevmode\hbox to3em{\hrulefill}\thinspace}

\end{document}